\documentclass[11pt, twoside]{article}
\usepackage{latexsym}
\usepackage{amsmath}
\usepackage{amssymb}
\usepackage[all]{xy}
\usepackage{amsfonts}
\usepackage{verbatim}
\usepackage{amsthm}
\usepackage{mathrsfs}\usepackage{array}
\usepackage{epsfig}
\usepackage{xy}
\usepackage{tensor}
\usepackage{array}
\usepackage{stmaryrd}
\usepackage{graphicx,color}
\usepackage{xcolor}
\usepackage[colorlinks=true,linkcolor=blue,citecolor=blue]{hyperref}
\usepackage{tikz}
\usetikzlibrary{arrows,calc}
\usepackage{etex}
\usepackage{mathdots}
\usepackage{float}
\usepackage{graphics}
\usepackage{pdflscape}
\usepackage{extarrows}
\usepackage{anysize,hyperref}
\input xypic
\xyoption{all}
\usepackage{bm}
\usepackage[perpage,symbol]{footmisc}
\usepackage{setspace}
\SelectTips{cm}{10}

\def\A{\mathscr{A}}

\def\C{\mathscr{C}}

\def\E{\mathbb{E}}
\def\s{\mathfrak{s}}
\def\Id{\mathrm{Id}}

\def\del{\delta}
\def\dr{\ar@{->}[r]}
\def\Im{\mbox{\rm Im}\,}\def\Ker{\mbox{\rm Ker}\,}

\begin{document}
\baselineskip=15pt
\title{\Large{\bf  Balanced pairs and recollements in extriangulated categories with negative first extensions\footnotetext{ Panyue Zhou was supported by the National Natural Science Foundation of China (Grant No. 11901190) and the Scientific Research Fund of Hunan Provincial Education Department (Grant No. 19B239).} }}
\medskip
\author{Jian He and Panyue Zhou}

\date{}

\maketitle
\def\blue{\color{blue}}
\def\red{\color{red}}

\newtheorem{theorem}{Theorem}[section]
\newtheorem{lemma}[theorem]{Lemma}
\newtheorem{corollary}[theorem]{Corollary}
\newtheorem{proposition}[theorem]{Proposition}
\newtheorem{conjecture}{Conjecture}
\theoremstyle{definition}
\newtheorem{definition}[theorem]{Definition}
\newtheorem{question}[theorem]{Question}
\newtheorem{notation}[theorem]{Notation}
\newtheorem{remark}[theorem]{Remark}
\newtheorem{remark*}[]{Remark}
\newtheorem{example}[theorem]{Example}
\newtheorem{example*}[]{Example}

\newtheorem{construction}[theorem]{Construction}
\newtheorem{construction*}[]{Construction}

\newtheorem{assumption}[theorem]{Assumption}
\newtheorem{assumption*}[]{Assumption}

\baselineskip=17pt
\parindent=0.5cm

\begin{abstract}
\baselineskip=16pt A notion of balanced pairs in an extriangulated category with a negative first extension is defined in this article. We prove that there exists a bijective correspondence between balanced pairs and proper classes $\xi$ with enough {$\xi$-projectives} and enough {$\xi$-injectives}. It can be regarded as a simultaneous generalization of Fu-Hu-Zhang-Zhu and Wang-Li-Huang. Besides, we show that if $(\mathcal A ,\mathcal B,\mathcal C)$ is a recollement of extriangulated categories, then
balanced pairs in $\mathcal B$ can induce balanced pairs in $\mathcal A$ and $\mathcal C$ under natural assumptions. As a application, this result gengralizes a result by Fu-Hu-Yao in abelian categories. Moreover, it highlights a new phenomena when it applied to triangulated categories.\\[2mm]
\textbf{Keywords:}  extriangulated categories; proper class; balanced pair; recollement.\\[2mm]
\textbf{ 2020 Mathematics Subject Classification:} 18G80; 18E10; 18E40.
\end{abstract}
\medskip

\pagestyle{myheadings}
\markboth{\rightline {\scriptsize J. He and P. Zhou\hspace{2mm}}}
         {\leftline{\scriptsize Balanced pairs and recollements in extriangulated categories with negative first extensions}}

\section{Introduction}
The recollement of triangulated categories was introduced first by Beilinson, Bernstein, and Deligne, see \cite{BBD}. A fundamental example of a recollement situation of abelian categories appeared in the construction of perverse sheaves by MacPherson and Vilonen \cite{MV}. Recollements of abelian and triangulated categories play an important role in
ring theory, representation theory and geometry of singular spaces.

Relative homological algebra has been formulated by Hochschild in categories of modules and later by Heller, Butler and Horrocks in more general categories with a relative abelian structure.
Beligiannis \cite{B} studied a homological algebra in a triangulated category  which parallels the homological algebra in an exact category in the sense of Quillen, by specifying a class of
triangles $\xi$ which is called a proper class of triangles.

The notion of extriangulated categories was introduced by Nakaoka and Palu in \cite{NP} as a simultaneous generalization of
exact categories and triangulated categories. Exact categories (abelian categories are also exact categories) and extension closed subcategories of an
extriangulated category are extriangulated categories, while there are some other examples of extriangulated categories which are neither exact nor triangulated, see \cite{NP,ZZ,HZZ,NP1}. Hence many results hold on exact categories and triangulated categories can be unified in the same framework. Recently, Adachi, Enomoto and Tsukamoto \cite{AET} introduced the notion of extriangulated categories with negative first extensions. They also showed that exact categories and triangulated categories naturally admit negative first extension structures.
Based on Beligiannis's idea,
Hu, Zhang and Zhou \cite{HZZ} defined a proper class of extriangulated categories, they proved that if an extriangulated category $(\mathcal{C},\E,\s)$ was equipped with a proper class of $\mathbb{E}$-triangles $\xi$, then $\mathcal{C}$ had a new extriangulated structure. This construction gives extriangulated categories which are neither exact nor triangulated.
Wang, Wei, and Zhang \cite{WWZ} introduced the recollement of extriangulated categories, which is a simultaneous generalization of recollements of abelian categories and triangulated categories.

Chen \cite{C}  introduced the notion of balanced pair of additive subcategories
in an abelian category. He studied relative homology with respect to balanced pairs in abelian category. Let $(\mathcal A ,\mathcal B,\mathcal C)$ be a recollement of abelian categories. Fu, Hu and Yao \cite{FHY} proved that balanced pairs in $\mathcal B$ can induce balanced pairs in $\mathcal A$ and $\mathcal C$. Wang, Li, and Huang \cite{WLH} showed that there exists a one-to-one correspondence between balanced pairs and Quillen exact structures $\xi$ with enough {$\xi$-projectives} and enough {$\xi$-injectives} for abelian categories.
Recently, the notion of balanced pair in triangulated categories was introduced by Fu, Hu, Zhang and Zhu \cite{FHZZ}, and proved a similar result to Wang, Li, and Huang \cite{WLH}. More precisely, they showed there exists a bijective correspondence between balanced pairs and proper classes $\xi$ with enough {$\xi$-projectives} and enough {$\xi$-injectives}. Inspired by this, we have a natural question whether their results of \cite{WLH} and \cite{FHZZ} can be unified under the framework of extriangulated categories. In this article, we give an affirmative answer in extriangulated categories with negative first extensions.

Let $\mathcal{C}$ be an extriangulated category with a negative first extension. Our first main result constructs a bijective correspondence between balanced pairs in $\mathcal C$ and proper classes $\xi$ with enough {$\xi$-projectives} and enough {$\xi$-injectives}, see Theorem \ref{main}. This unifies their results of  Fu, Hu, Zhang and Zhu \cite{FHZZ} and Wang, Li and Huang \cite{WLH} in the framework of extriangulated categories with negative first extensions. Suppose that $\mathcal B$ admits a recollement relative to extriangulated categories $\mathcal A$ and $\mathcal C$. Our second main result show that balanced pairs in $\mathcal B$ can induce balanced pairs in $\mathcal A$ and $\mathcal C$ under natural assumptions, see Theorem \ref{www}. This generalizes a result of Fu, Hu and Yao \cite{FHY}.

This article is organized as follows. In Section 2, we give some terminologies
and some preliminary results. In Section 3, we prove our first main
results. In Section 4, we prove our second main
result.

\section{Preliminaries}
We briefly recall some definitions and basic properties of extriangulated categories from \cite{NP}.
We omit some details here, but the reader can find them in \cite{NP}.

Let $\mathcal{C}$ be an additive category equipped with an additive bifunctor
$$\mathbb{E}: \mathcal{C}^{\rm op}\times \mathcal{C}\rightarrow {\rm Ab},$$
where ${\rm Ab}$ is the category of abelian groups. For any objects $A, C\in\mathcal{C}$, an element $\delta\in \mathbb{E}(C,A)$ is called an $\mathbb{E}$-extension.
Let $\mathfrak{s}$ be a correspondence which associates an equivalence class $$\mathfrak{s}(\delta)=\xymatrix@C=0.8cm{[A\ar[r]^x
 &B\ar[r]^y&C]}$$ to any $\mathbb{E}$-extension $\delta\in\mathbb{E}(C, A)$. This $\mathfrak{s}$ is called a {\it realization} of $\mathbb{E}$, if it makes the diagrams in \cite[Definition 2.9]{NP} commutative.
 A triplet $(\mathcal{C}, \mathbb{E}, \mathfrak{s})$ is called an {\it extriangulated category} if it satisfies the following conditions.
\begin{itemize}
\item $\mathbb{E}\colon\mathcal{C}^{\rm op}\times \mathcal{C}\rightarrow \rm{Ab}$ is an additive bifunctor.

\item $\mathfrak{s}$ is an additive realization of $\mathbb{E}$.

\item $\mathbb{E}$ and $\mathfrak{s}$  satisfy the compatibility conditions in \cite[Definition 2.12]{NP}.
\end{itemize}

We collect the following terminology from \cite{NP}.

\begin{definition}
Let $(\mathcal{C},\E,\s)$ be an extriangulated category.
\begin{itemize}
\item[(1)] A sequence $A\xrightarrow{~x~}B\xrightarrow{~y~}C$ is called a {\it conflation} if it realizes some $\E$-extension $\del\in\E(C,A)$.
    In this case, $x$ is called an {\it inflation} and $y$ is called a {\it deflation}.

\item[(2)] If a conflation  $A\xrightarrow{~x~}B\xrightarrow{~y~}C$ realizes $\delta\in\mathbb{E}(C,A)$, we call the pair $( A\xrightarrow{~x~}B\xrightarrow{~y~}C,\delta)$ an {\it $\E$-triangle}, and write it in the following way.
$$A\overset{x}{\longrightarrow}B\overset{y}{\longrightarrow}C\overset{\delta}{\dashrightarrow}$$
We usually do not write this $``\delta"$ if it is not used in the argument.

\item[(3)] Let $A\overset{x}{\longrightarrow}B\overset{y}{\longrightarrow}C\overset{\delta}{\dashrightarrow}$ and $A^{\prime}\overset{x^{\prime}}{\longrightarrow}B^{\prime}\overset{y^{\prime}}{\longrightarrow}C^{\prime}\overset{\delta^{\prime}}{\dashrightarrow}$ be any pair of $\E$-triangles. If a triplet $(a,b,c)$ realizes $(a,c)\colon\delta\to\delta^{\prime}$, then we write it as
$$\xymatrix{
A \ar[r]^x \ar[d]^a & B\ar[r]^y \ar[d]^{b} & C\ar@{-->}[r]^{\del}\ar[d]^c&\\
A'\ar[r]^{x'} & B' \ar[r]^{y'} & C'\ar@{-->}[r]^{\del'} &}$$
and call $(a,b,c)$ a {\it morphism of $\E$-triangles}.

\end{itemize}
\end{definition}
\begin{proposition}\label{exact} \rm{\cite[Proposition 3.3]{NP}}
Let $\mathcal{C}$ be an extriangulated category. For any $\mathbb{E}$-triangle $A\stackrel{}{\longrightarrow}B\stackrel{}{\longrightarrow}C\stackrel{\delta}\dashrightarrow$, we have the following exact sequences:
$$\mathcal{C}(C,-)\rightarrow \mathcal{C}(B,-)\rightarrow \mathcal{C}(A,-)\stackrel{~}\rightarrow\mathbb{E}(C,-)\rightarrow\mathbb{E}(B,-);$$
$$\mathcal{C}(-,A)\rightarrow \mathcal{C}(-,B)\rightarrow \mathcal{C}(-,C)\stackrel{~}\rightarrow\mathbb{E}(-,A)\rightarrow\mathbb{E}(-,B).$$
\end{proposition}
Let us recall the definition of a negative first extension structure on an extriangulated category form \cite{AET}.
\begin{definition}\label{qq}
Let $\mathcal{C}$ be an extriangulated category. A \emph{negative first extension structure} on $\mathcal{C}$ consists of the following data:
\begin{itemize}
\item[(NE1)] $\mathbb{E}^{-1}: \mathcal{C}^{\rm op}\times \mathcal{C}\rightarrow {\rm Ab}$ is an additive bifunctor.
\item[(NE2)] For each $\delta\in \mathbb{E}(C,A)$, there exist two natural transformations
\begin{align}
&\delta_{\sharp}^{-1}: \mathbb{E}^{-1}(-,C)\to \mathcal{C}(-,A),\notag\\
&\delta^{\sharp}_{-1}: \mathbb{E}^{-1}(A,-)\to \mathcal{C}(C,-)\notag
\end{align}
such that for each $\mathbb{E}$-triangle $A\stackrel{f}{\longrightarrow}B\stackrel{g}{\longrightarrow}C\stackrel{\delta}\dashrightarrow$

$\xymatrix@C=1cm{\mathbb{E}^{-1}(C, -)\ar[r]^{\mathbb{E}^{-1}(g, -)}&\mathbb{E}^{-1}(B, -)\ar[r]^{\mathbb{E}^{-1}(f, -)}&\mathbb{E}^{-1}(A, -)\ar[r]^{\delta^{\sharp}_{-1}}&\mathcal{C}(C, -)\ar[r]^{\mathcal{C}(g, -)}&\mathcal{C}(B, -);}$

$\xymatrix@C=1cm{\mathbb{E}^{-1}(-, A)\ar[r]^{\mathbb{E}^{-1}(-, f)}&\mathbb{E}^{-1}(-, B)\ar[r]^{\mathbb{E}^{-1}(-, g)}&\mathbb{E}^{-1}(-, C)\ar[r]^{\delta_{\sharp}^{-1}}&\mathcal{C}(-, A)\ar[r]^{\mathcal{C}(-, f)}&\mathcal{C}(-, B);}$

are exact.
\end{itemize}
Then we call $\mathcal{C}=(\mathcal{C},\mathbb{E},\mathfrak{s},\mathbb{E}^{-1})$ an \emph{extriangulated category with a negative first extension}.
\end{definition}

\begin{remark}\label{111}
Typical examples of extriangulated categories are triangulated categories and exact categories (see \cite[Example 2.13]{NP}).
Adachi, Enomoto and Tsukamoto \cite{AET} show that both categories naturally admit negative first extension structures (see \cite[Example 2.4]{AET}).
\end{remark}

Let $(\mathcal C,\E,\s)$ be an extrangulated category. Since $\E$ is a bifunctor, for any $a\in\mathcal C(A,A')$ and $c\in\mathcal C(C',C)$, we have $\E$-extensions
$$ \E(C,a)(\del)\in\E(C,A')\ \ \text{and}\ \ \ \E(c,A)(\del)\in\E(C',A). $$
We simply denote them by $a_\ast\delta$ and $c^\ast\delta$.

\begin{lemma}\label{lem1} \emph{\cite[Proposition 3.15]{NP}} Let $(\mathcal{C}, \mathbb{E},\mathfrak{s})$ be an extriangulated category, $$\xymatrix@C=2em{A_1\ar[r]^{x_1}&B_1\ar[r]^{y_1}&C\ar@{-->}[r]^{\delta_1}&}~\mbox{and}~ \xymatrix@C=2em{A_2\ar[r]^{x_2}&B_2\ar[r]^{y_2}&C\ar@{-->}[r]^{\delta_2}&}$$ be any pair of $\mathbb{E}$-triangles. Then there exists a commutative diagram
in $\mathcal{C}$
$$\xymatrix{
    & A_2\ar[d]_{m_2} \ar@{=}[r] & A_2 \ar[d]^{x_2} \\
  A_1 \ar@{=}[d] \ar[r]^{m_1} & M \ar[d]_{e_2} \ar[r]^{e_1} & B_2\ar[d]^{y_2} \\
  A_1 \ar[r]^{x_1} & B_1\ar[r]^{y_1} & C   }
  $$
  which satisfies $\mathfrak{s}(y^*_2\delta_1)=\xymatrix@C=2em{[A_1\ar[r]^{m_1}&M\ar[r]^{e_1}&B_2]}$ and
  $\mathfrak{s}(y^*_1\delta_2)=\xymatrix@C=2em{[A_2\ar[r]^{m_2}&M\ar[r]^{e_2}&B_1].}$
\end{lemma}

We now review the notion of proper classes of $\mathbb{E}$-triangles  and its related properties from \cite{HZZ}. From now on, assume that $(\mathcal C,\E,\s)$ is an extrangulated category.

A class of $\mathbb{E}$-triangles $\xi$ is {\it closed under base change} if for any $\mathbb{E}$-triangle $$\xymatrix@C=2em{A\ar[r]^x&B\ar[r]^y&C\ar@{-->}[r]^{\delta}&\in\xi}$$ and any morphism $c\colon C' \to C$, then any $\mathbb{E}$-triangle  $\xymatrix@C=2em{A\ar[r]^{x'}&B'\ar[r]^{y'}&C'\ar@{-->}[r]^{c^*\delta}&}$ belongs to $\xi$.

Dually, a class of  $\mathbb{E}$-triangles $\xi$ is {\it closed under cobase change} if for any $\mathbb{E}$-triangle $$\xymatrix@C=2em{A\ar[r]^x&B\ar[r]^y&C\ar@{-->}[r]^{\delta}&\in\xi}$$ and any morphism $a\colon A \to A'$, then any $\mathbb{E}$-triangle  $\xymatrix@C=2em{A'\ar[r]^{x'}&B'\ar[r]^{y'}&C\ar@{-->}[r]^{a_*\delta}&}$ belongs to $\xi$.

A class of $\mathbb{E}$-triangles $\xi$ is called {\it saturated} if in the situation of Lemma \ref{lem1}, whenever {
 $\xymatrix@C=2em{A_2\ar[r]^{x_2}&B_2\ar[r]^{y_2}&C\ar@{-->}[r]^{\delta_2 }&}$
 and $\xymatrix@C=2em{A_1\ar[r]^{m_1}&M\ar[r]^{e_1}&B_2\ar@{-->}[r]^{y_2^{\ast}\delta_1}&}$ }
 belong to $\xi$, then the  $\mathbb{E}$-triangle $\xymatrix@C=2em{A_1\ar[r]^{x_1}&B_1\ar[r]^{y_1}&C\ar@{-->}[r]^{\delta_1 }&}$  belongs to $\xi$.

An $\mathbb{E}$-triangle $\xymatrix@C=2em{A\ar[r]^x&B\ar[r]^y&C\ar@{-->}[r]^{\delta}&}$ is called {\it split} if $\delta=0$. It is easy to see that it is split if and only if $x$ is section or $y$ is retraction. The full subcategory  consisting of the split $\mathbb{E}$-triangles will be denoted by $\Delta_0$.

  \begin{definition}{\rm \cite[Definition 3.1]{HZZ}} \label{def:proper class}{\rm  Let $\xi$ be a class of $\mathbb{E}$-triangles which is closed under isomorphisms. $\xi$ is called a {\it proper class} of $\mathbb{E}$-triangles if the following conditions hold:

  (1) $\xi$ is closed under finite coproducts and $\Delta_0\subseteq \xi$.

  (2) $\xi$ is closed under base change and cobase change.

  (3) $\xi$ is saturated.}

  \end{definition}
We fix a proper class $\xi$ of $\mathbb{E}$-triangles in the $\mathcal{C}$.

 \begin{definition}{\rm \cite[Definition 4.1]{HZZ}}  {\rm An object $P\in\mathcal{C}$  is called {\it $\xi$-projective}  if for any $\mathbb{E}$-triangle $$\xymatrix{A\ar[r]^x& B\ar[r]^y& C \ar@{-->}[r]^{\delta}& }$$ in $\xi$, the induced sequence of abelian groups $\xymatrix@C=0.6cm{0\ar[r]& \mathcal{C}(P,A)\ar[r]& \mathcal{C}(P,B)\ar[r]&\mathcal{C}(P,C)\ar[r]& 0}$ is exact. Dually, we have the definition of {\it $\xi$-injective}.}
\end{definition}
We denote by $\mathcal{P(\xi)}$ (respectively $\mathcal{I(\xi)}$) the class of $\xi$-projective (respectively $\xi$-injective) objects of $\mathcal{C}$. It follows from the definition that this subcategory $\mathcal{P}(\xi)$ and $\mathcal{I}(\xi)$ are full, additive, closed under isomorphisms and direct summands.

 An extriangulated  category $(\mathcal{C}, \mathbb{E}, \mathfrak{s})$  is said to  have {\it  enough
$\xi$-projectives} \ (respectively {\it  enough $\xi$-injectives}) provided that for each object $A$ there exists an $\mathbb{E}$-triangle $\xymatrix@C=2.1em{K\ar[r]& P\ar[r]&A\ar@{-->}[r]& }$ (respectively $\xymatrix@C=2em{A\ar[r]& I\ar[r]& K\ar@{-->}[r]&}$) in $\xi$ with $P\in\mathcal{P}(\xi)$ (respectively $I\in\mathcal{I}(\xi)$).

\begin{definition}{\rm \cite[Definition 4.4]{HZZ}}  {\rm An unbounded complex $\mathbf{X}$ is called {\it $\xi$-exact} complex  if $\mathbf{X}$ is a diagram $$\xymatrix@C=2em{\cdots\ar[r]&X_1\ar[r]^{d_1}&X_0\ar[r]^{d_0}&X_{-1}\ar[r]&\cdots}$$ in $\mathcal{C}$ such that for each integer $n$, there exists an $\mathbb{E}$-triangle $\xymatrix@C=2em{K_{n+1}\ar[r]^{g_n}&X_n\ar[r]^{f_n}&K_n\ar@{-->}[r]^{\delta_n}&}$ in $\xi$ and $d_n=g_{n-1}f_n$.}
\end{definition}
\begin{definition}{\rm \cite[Definition 4.5]{HZZ}}  {\rm Let $\mathcal{W}$ be a class of objects in $\mathcal{C}$. An $\mathbb{E}$-triangle $$\xymatrix@C=2em{A\ar[r]& B\ar[r]& C\ar@{-->}[r]& }$$ in $\xi$ is called to be
{\it $\mathcal{C}(-,\mathcal{W})$-exact} (respectively
{\it $\mathcal{C}(\mathcal{W},-)$-exact}) if for any $W\in\mathcal{W}$, the induced sequence of abelian group $\xymatrix@C=2em{0\ar[r]&\mathcal{C}(C,W)\ar[r]&\mathcal{C}(B,W)\ar[r]&\mathcal{C}(A,W)\ar[r]& 0}$ (respectively \\ $\xymatrix@C=2em{0\ar[r]&\mathcal{C}(W,A)\ar[r]&\mathcal{C}(W,B)\ar[r]&\mathcal{C}(W,C)\ar[r]&0}$) is exact in ${\rm Ab}$}.
\end{definition}

\section{On the relation between proper classes of $\mathbb{E}$-triangles and balanced pairs }
In the rest of this article, unless otherwise stated, we always regard $\mathcal{C}$ as an extriangulated category with a negative first extension.
\begin{definition}\label{dd1}
Let $\mathcal{C}$ be an  extriangulated category. A subcategory $\mathcal{T}$ of $\mathcal{C}$ is called
\emph{strongly exact-contravariantly finite}, if for any object $C\in\mathcal{C}$, there exists an $\E$-triangle
$$\xymatrix{K\ar[r]&T\ar[r]^{g}&C\ar@{-->}[r]^{\del}&,}~~~~~T\in\mathcal{T},$$
such that the induced sequence of abelian group $$0\rightarrow\mathcal{C}(T',K)\rightarrow\mathcal{C}(T',T)\rightarrow\mathcal{C}(T',C)\rightarrow 0$$ is exact in ${\rm Ab},$ where $T'\in\mathcal{T}$. i.e. the above $\E$-triangle is {\it $\mathcal{C}(\mathcal{T},-)$-exact}.

Dually, a subcategory $\mathcal{T}$ of $\mathcal{C}$ is called
\emph{strongly exact-covariantly  finite}, if for any object $C\in\mathcal{C}$, there exists an $\E$-triangle
$$\xymatrix{C\ar[r]^{f}&T\ar[r]&L\ar@{-->}[r]^{\del'}&,}~~~~~T\in\mathcal{T},$$
such that the induced sequence of abelian group $$0\rightarrow\mathcal{C}(L,T')\rightarrow\mathcal{C}(T,T')\rightarrow\mathcal{C}(C,T')\rightarrow 0$$ is exact in ${\rm Ab},$ where $T'\in\mathcal{T}$. i.e. the above $\E$-triangle is {\it $\mathcal{C}(-,\mathcal{T})$-exact}.

\end{definition}

\begin{remark}\label{dd1}
 Let $\mathcal{T}\subseteq\mathcal{C}$ be a $\Sigma$-stable $(i.e. \Sigma\mathcal{T}=\mathcal{T})$ subcategory of triangulated category $\mathcal{C}$. $\mathcal{T}$ is \emph{strongly exact-contravariantly  finite} in $\mathcal{C}$ if and only if $\mathcal{T}$ is \emph{contravariantly  finite} in $\mathcal{C}$, $\mathcal{T}$ is \emph{strongly exact-covariantly  finite} in $\mathcal{C}$ if and only if $\mathcal{T}$ is \emph{covariantly  finite} in $\mathcal{C}$ by Lemma 2.6 in \cite{HYF}.
\end{remark}
Assume that $\mathcal{X}$ is a full additive subcategory of $\mathcal{C}$. For any $\E$-triangle in $\mathcal{C}$
$$T:\xymatrix{X\ar[r]^{f}&Y\ar[r]^{g}&Z\ar@{-->}[r]^{\del}&,}$$
We write \begin{align*}
\xi_{\mathcal{X}}=\left\{
 \begin{array}{c|c}
 \textrm{the classes of } &   \textrm{the }\E\textrm{-triangle }T\textrm{ is}   \\
 \E\textrm{-triangle }T  &\mathcal C(\mathcal{X},-)\textrm{-exact}
 \end{array}
\right\}
\end{align*}

\begin{proposition}\label{gg}
 Let $\mathcal{X}$ be a full additive subcategory of $\mathcal{C}$ which is closed under direct summands. Then $\xi_{\mathcal{X}}$ is a proper class of $\mathcal{C}$. Moreover, $\mathcal{X}$ is \emph{strongly exact-contravariantly finite} in $\mathcal{C}$ if and only if $\mathcal{C}$  has {\it  enough
$\xi_{\mathcal{X}}$-projectives} and ${\mathcal{X}=\mathcal{P}(\xi_{\mathcal{X}}})$.
\begin{proof}It is easy to see that $\xi_{\mathcal{X}}$ is closed under isomorphisms, finite coproducts and containing all split $\E$-triangles. Next we claim that $\xi_{\mathcal{X}}$ is closed under base change and cobase change. Consider the following commutative diagram of $\E$-triangles:
$$\xymatrix{
A \ar[r]^{x'} \ar@{=}[d] & B'\ar[r]^{y'} \ar[d]^{b} & C'\ar@{-->}[r]^{c^*\del}\ar[d]^c&\\
A\ar[r]^{x} & B \ar[r]^{y} & C\ar@{-->}[r]^{\del} &,}$$
then we have the the following commutative by Proposition \ref{exact}
$$\xymatrix@C=1.2cm{
\mathcal{C}(X,A)\ar[r]^{\mathcal{C}(X, x')}\ar@{=}[d] &\mathcal{C}(X,B')\ar[r]^{\mathcal{C}(X, y')}\ar[d]^{\mathcal{C}(X,b)} & \mathcal{C}(X,C') \ar[r] \ar[d]^{\mathcal{C}(X,c)}& \mathbb{E}(X,A)\ar[r]^{\mathbb{E}(X, x')}\ar@{=}[d]& \mathbb{E}(X,B')\ar[d]^{\mathbb{E}(X,b)}\\
\mathcal{C}(X,A)\ar[r]^{\mathcal{C}(X, x)} &\mathcal{C}(X,B)\ar[r]^{\mathcal{C}(X, y)}\ar[r] &  \mathcal{C}(X,C) \ar[r]& \mathbb{E}(X,A)\ar[r]^{\mathbb{E}(X, x)}& \mathbb{E}(X,B)
}$$
for any $X\in{\mathcal{X}}$. If the $\E$-triangle $\xymatrix{A\ar[r]^{x}&B\ar[r]^{y}&C\ar@{-->}[r]^{\del}&}$ belongs to $\xi_{\mathcal{X}}$, then $\mathcal{C}(X,b)\mathcal{C}(X, x')=\mathcal{C}(X, x)$  is monic, thus $\mathcal{C}(X, x')$ is monic. Similarly, since $\mathcal{C}(X, y)$ is epic, then $\mathbb{E}(X,b)\mathbb{E}(X, x')=\mathbb{E}(X, x)$ is monic, thus $\mathbb{E}(X, x')$ is monic, i.e. $\mathcal{C}(X, y')$ is epic. Hence $\xi_{\mathcal{X}}$ is closed under base change.

Consider the following commutative diagram of $\E$-triangles:
$$\xymatrix{
A \ar[r]^{x} \ar[d]^{a} & B\ar[r]^{y} \ar[d]^{b} & C\ar@{-->}[r]^{\del}\ar@{=}[d]&\\
A'\ar[r]^{x'} & B' \ar[r]^{y'} & C\ar@{-->}[r]^{c_*\del} &,}$$
then we have the the following commutative by (NE2)
$$\xymatrix@C=1.2cm{
\mathbb{E}^{-1}(X,B)\ar[r]^{\mathbb{E}^{-1}(X, x)}\ar[d]^{\mathbb{E}^{-1}(X,b)} &\mathbb{E}^{-1}(X,C)\ar[r]^{}\ar@{=}[d]& \mathcal{C}(X,A) \ar[r]^{\mathcal{C}(X, x)} \ar[d]^{\mathcal{C}(X,a)}& \mathcal{C}(X,B)\ar[r]^{\mathcal{C}(X, y)}\ar[d]^{\mathcal{C}(X,b)}& \mathcal{C}(X,C)\ar@{=}[d]\\
\mathbb{E}^{-1}(X,B')\ar[r]^{\mathbb{E}^{-1}(X, x')} &\mathbb{E}^{-1}(X,C)\ar[r]&  \mathcal{C}(X,A') \ar[r]^{\mathcal{C}(X, x')}& \mathcal{C}(X,B')\ar[r]^{\mathcal{C}(X, y')}& \mathcal{C}(X,C)
}$$
for any $X\in{\mathcal{X}}$. If the $\E$-triangle $\xymatrix{A\ar[r]^{x}&B\ar[r]^{y}&C\ar@{-->}[r]^{\del}&}$ belongs to $\xi_{\mathcal{X}}$, then $\mathcal{C}(X,y')\mathcal{C}(X, b)=\mathcal{C}(X, y)$  is epic, thus $\mathcal{C}(X,y')$ is epic. Similarly, since $\mathcal{C}(X, x)$ is monic, then $\mathbb{E}^{-1}(X,x')\mathbb{E}^{-1}(X, b)=\mathbb{E}^{-1}(X, x)$ is epic, thus $\mathbb{E}^{-1}(X, x')$ is epic, i.e. $\mathcal{C}(X, x')$ is monic. Hence $\xi_{\mathcal{X}}$ is closed under cobase change.

Consider the following commutative diagram of $\E$-triangles of Lemma \ref{lem1}:
$$\xymatrix{
    & A_2\ar[d]_{m_2} \ar@{=}[r] & A_2 \ar[d]^{x_2} \\
  A_1 \ar@{=}[d] \ar[r]^{m_1} & M \ar[d]_{e_2} \ar[r]^{e_1} & B_2\ar[d]^{y_2} \\
  A_1 \ar[r]^{x_1} & B_1\ar[r]^{y_1} & C   }
  $$
then we have the the following commutative by (NE2)
$$\xymatrix@C=1.2cm{
\mathbb{E}^{-1}(X,M)\ar[r]^{\mathbb{E}^{-1}(X, e_1)}\ar[d]^{\mathbb{E}^{-1}(X,e_2)} &\mathbb{E}^{-1}(X,B_2)\ar[r]^{}\ar[d]^{\mathbb{E}^{-1}(X, y_2)}& \mathcal{C}(X,A_1) \ar[r]^{\mathcal{C}(X, m_1)} \ar@{=}[d]& \mathcal{C}(X,M) \ar[r]^{\mathcal{C}(X, e_1)} \ar[d]^{\mathcal{C}(X,e_2)}& \mathcal{C}(X,B_2)\ar[d]^{\mathcal{C}(X,y_2)}\\
\mathbb{E}^{-1}(X,B_1)\ar[r]^{\mathbb{E}^{-1}(X, y_1)} &\mathbb{E}^{-1}(X,C)\ar[r]&  \mathcal{C}(X,A_1) \ar[r]^{\mathcal{C}(X, x_1)}&  \mathcal{C}(X,B_1) \ar[r]^{\mathcal{C}(X, y_1)}& \mathcal{C}(X,C)&
}$$
for any $X\in{\mathcal{X}}$. If the $\E$-triangle $\xymatrix{A_1\ar[r]^{m_1}&M\ar[r]^{e_1}&B_2\ar@{-->}[r]^{}&}$ and $\xymatrix{A_2\ar[r]^{x_2}&B_2\ar[r]^{y_2}&C\ar@{-->}[r]^{}&}$  belong to $\xi_{\mathcal{X}}$, then we have $\mathcal{C}(X,e_1), \mathbb{E}^{-1}(X,e_1),\mathcal{C}(X,y_2)$ and $\mathbb{E}^{-1}(X,y_2)$ are epic. It is easy to see $\mathbb{E}^{-1}(X,y_1)$ and $\mathcal{C}(X,y_1)$ are epic, which implies that $\E$-triangle $\xymatrix{A_1\ar[r]^{x_1}&B_1\ar[r]^{y_1}&C\ar@{-->}[r]^{}&}$ belongs to $\xi_{\mathcal{X}}$, Hence $\xi_{\mathcal{X}}$ is saturated.

Finally, if $\mathcal{C}$  has {\it  enough $\xi_{\mathcal{X}}$-projectives} and ${\mathcal{X}=\mathcal{P}(\xi_{\mathcal{X}}})$ then it is easy to check that $\mathcal{X}$ is \emph{strongly exact-contravariantly finite} in $\mathcal{C}$. Conversely, assume that $\mathcal{X}$ is \emph{strongly exact-contravariantly finite} in $\mathcal{C}$. It is clear that ${\mathcal{X}\subseteq\mathcal{P}(\xi_{\mathcal{X}}})$. On the other hand, let ${P\in\mathcal{P}(\xi_{\mathcal{X}}})$. Then there exists an $\E$-triangle $\xymatrix{K\ar[r]^{}&X\ar[r]^{}&P\ar@{-->}[r]^{}&}$ in ${\xi_{\mathcal{X}}}$ with $X\in{{\mathcal{X}}}$. Thus this $\E$-triangle is split, and hence $P\in{{\mathcal{X}}}$. This completes the proof.
\end{proof}
\end{proposition}

Assume that $\mathcal{Y}\subseteq\mathcal{C}$ is a full additive subcategory of $\mathcal{C}$. For any $\E$-triangle in $\mathcal{C}$ $$T:\xymatrix{X\ar[r]^{f}&Y\ar[r]^{g}&Z\ar@{-->}[r]^{\del}&,}$$
We write
\begin{align*}
\xi^{\mathcal{Y}}=\left\{
 \begin{array}{c|c}
 \textrm{the classes of } &   \textrm{the }\E\textrm{-triangle }T\textrm{ is}   \\
 \E\textrm{-triangle }T  & C(-,\mathcal{Y})\textrm{-exact}
 \end{array}
\right\}
\end{align*}

\begin{proposition}\label{gbg}
 Let $\mathcal{Y}$ be a full additive subcategory of $\mathcal{C}$ which is closed under direct summands. Then $\xi^{\mathcal{Y}}$ is a proper class of $\mathcal{C}$. Moreover, $\mathcal{Y}$ is \emph{strongly exact-covariantly finite} in $\mathcal{C}$ if and only if $\mathcal{C}$ has {\it  enough
$\xi^{\mathcal{Y}}$-injectives} and ${\mathcal{Y}=\mathcal{I}(\xi^{\mathcal{Y}}})$.
\begin{proof}
This follows from the dual of Proposition \ref{gg} and Lemma 4.10 (1) in \rm{\cite{HZZ}}.
\end{proof}

\end{proposition}
\begin{definition}\label{77}
Let $C\in\mathcal{C}$,  $\mathcal{X}$ and $\mathcal{Y}$ be full additive subcategories of $\mathcal{C}$.
\begin{itemize}
 \rm\item[(1)] An $\mathcal{X}$-resolution of $C$ is a diagram $X^\bullet\rightarrow C$ such that $ \xymatrix{X^\bullet:= \cdots\ar[r]&X_1\ar[r]^{}&X_0\ar[r]^{}\ar[r]&0}$ is a complex with $X_i\in\mathcal{X}$ for all $i\geq0$, and $ \xymatrix{\cdots\ar[r]&X_1\ar[r]^{}&X_0\ar[r]^{}\ar[r]&C}$ is a $\xi_{\mathcal{X}}$-exact complex. Moreover, the $\mathcal{X}$-resolution $X^\bullet\rightarrow C$ of  $C$ is called $C(-,\mathcal{Y})$-exact if its $\xi_{\mathcal{X}}$-exact complex is $C(-,\mathcal{Y})$-exact.
 \rm\item[(2)] An $\mathcal{Y}$-coresolution of $C$ is a diagram $ C\rightarrow Y^\bullet$ such that $ \xymatrix@C=0.5cm{Y^\bullet:= 0\ar[r]&Y_0\ar[r]^{}&Y_{-1}\ar[r]^{}\ar[r]&\cdots}$ is a complex with $Y_i\in\mathcal{Y}$ for all $i\leq0$, and $ \xymatrix{C\ar[r]&Y_0\ar[r]^{}&Y_1\ar[r]^{}\ar[r]&\cdots}$ is a $\xi^{\mathcal{Y}}$-exact complex. Moreover, the $\mathcal{Y}$-resolution $ C\rightarrow Y^\bullet$ of $C$ is called $C(\mathcal{X},-)$-exact if its $\xi^{\mathcal{Y}}$-exact complex is  $C(\mathcal{X},-)$-exact.

 \end{itemize}

\end{definition}
\begin{remark}\label{xxx}
Using standard arguments from relative homological algebra, one can prove a version of the comparison theorem for $\mathcal{X}$-resolution (resp. $\mathcal{Y}$-coresolution). It follows that any two  $\mathcal{X}$-resolutions (resp. $\mathcal{Y}$-coresolutions) of an object $C$ are homotopy equivalences.
\end{remark}

\begin{proposition}\label{aag}
 Let $\mathcal{X}$ and $\mathcal{Y}$ be full additive subcategory of $\mathcal{C}$ which are closed under direct summands. Then the following statements are equivalent:
 \begin{itemize}
\rm\item[(1)] $\xi_{\mathcal{X}}=\xi^{\mathcal{Y}},{\mathcal{X}=\mathcal{P}(\xi_{\mathcal{X}}}),{\mathcal{Y}=\mathcal{I}(\xi^{\mathcal{X}}})$, and every object in $\mathcal{C}$ has enough {\it  enough $\xi_{\mathcal{X}}$-projectives} and enough {$\xi^{\mathcal{Y}}$-injectives}.
 \rm\item[(2)] The pair $(\mathcal{X},\mathcal{Y})$ satisfies
\begin{itemize}
  \item[(B1)] ${\mathcal{X}}$ is \emph{strongly exact-contravariantly finite} and $\mathcal{Y}$ is \emph{strongly exact-covariantly finite} in $\mathcal{C}$.

  \item[(B2)] For any object ${M\in\mathcal{C}}$, there exists an $\mathcal{X}$-resolution $X^{\bullet}\rightarrow{M}$ such that it is $C(-,\mathcal{Y})$-exact.

 \item[(B3)] For any object ${N\in\mathcal{C}}$, there exists a $\mathcal{Y}$-coresolution $N\rightarrow Y^{\bullet}$ such that it is $C(\mathcal{X},-)$-exact.
\end{itemize}
 \end{itemize}

\begin{proof}The proof is an adaption of \cite[Proposition 3.7]{FHZZ}.

$(1)\Rightarrow(2)$ It follows from Proposition \ref{gg} and Proposition \ref{gbg}.

$(2)\Rightarrow(1)$ By Proposition \ref{gg} and Proposition \ref{gbg}, it suffices to show $\xi_{\mathcal{X}}=\xi^{\mathcal{Y}}$.
Let $\xymatrix{A\ar[r]^{f}&B\ar[r]^{g}&C\ar@{-->}[r]^{}&}$ be an $\E$-triangle in ${\xi_{\mathcal{X}}}$. By hypothesis, there is an $\mathcal{X}$-resolution $X^{\bullet}\rightarrow{C}$ of ${C}$, such that it is $C(-,\mathcal{Y})$-exact. Then there exists a $\xi_{\mathcal{X}}$-exact complex
$$ \xymatrix{\cdots\ar[r]&X_1\ar[r]^{d_1}&X_0\ar[r]^{d_0}\ar[r]&C}$$ in ${\mathcal{C}}$ which is $C(-,\mathcal{Y})$-exact. This gives us an $\E$-triangle $\xymatrix{K_1\ar[r]^{g_0}&X_0\ar[r]^{f_0}&C\ar@{-->}[r]^{}&}$  which is $C(-,\mathcal{Y})$-exact, and hence we have the following commutative diagram of $\E$-triangles
$$\xymatrix{
K_1 \ar[r]^{g_0} \ar@{-->}[d]^{\alpha} & X_0\ar[r]^{f_0} \ar@{-->}[d]^{\beta} & C\ar@{-->}[r]^{}\ar@{=}[d]&\\
A\ar[r]^{f} & B \ar[r]^{g} & C\ar@{-->}[r]^{} &.}$$
Let $Y\in\mathcal{Y}$, applying $\mathcal C(-,\mathcal{Y})$ to the commutative diagram above, we have the following commutative diagram by Proposition \ref{exact}
$$\xymatrix@C=1.2cm{
\mathcal{C}(C,Y)\ar[r]^{\mathcal{C}(g, Y)}\ar@{=}[d] &\mathcal{C}(B,Y)\ar[r]^{\mathcal{C}(f, Y)}\ar[d]^{\mathcal{C}(\beta,Y)} & \mathcal{C}(A,Y) \ar[r] \ar[d]^{\mathcal{C}(\alpha,Y)}& \mathbb{E}(C,Y)\ar[r]^{\mathbb{E}(g, Y)}\ar@{=}[d]& \mathbb{E}(B,Y)\ar[d]^{\mathbb{E}(\beta,Y)}\\
\mathcal{C}(C,Y)\ar[r]^{\mathcal{C}(f_0, Y)} &\mathcal{C}(X_0,Y)\ar[r]^{\mathcal{C}(g_0, Y)}\ar[r] &  \mathcal{C}(K_1,Y) \ar[r]& \mathbb{E}(C,Y)\ar[r]^{\mathbb{E}(f_0, Y)}& \mathbb{E}(X_0,Y).
}$$
Note that ${\mathcal{C}(f_0, Y)}$ and ${\mathbb{E}(f_0, Y)}$ are monic, then we have that ${\mathcal{C}(g, Y)}$ and ${\mathbb{E}(g, Y)}$ are monic. Thus the $\E$-triangle $\xymatrix{A\ar[r]^{f}&B\ar[r]^{g}&C\ar@{-->}[r]^{}&}$ in ${\xi_{\mathcal{X}}}$ is $C(-,\mathcal{Y})$-exact and it belongs to $\xi^{\mathcal{Y}}$. This implies that $\xi_{\mathcal{X}}\subseteq\xi^{\mathcal{Y}}$. Dually, one can show that $\xi^{\mathcal{Y}}\subseteq\xi_{\mathcal{X}}$.
\end{proof}
\end{proposition}

Now we introduce the notion of balanced pairs.
\begin{definition}\label{bb} Suppose that $\mathcal{X}$ and $\mathcal{Y}$ are full additive subcategories of $\mathcal{C}$ which are closed under direct summands. The pair $(\mathcal{X},\mathcal{Y})$ is called a balanced pair if it satisfies one of the equivalence conditions of
of Proposition \ref{aag}.
\end{definition}

\begin{remark}\label{222}
(1) By Remark \ref{dd1}, if the categoriy $\mathcal{C}$ are triangulated, $\mathcal{X}$ and $\mathcal{Y}$ are $\Sigma$-stable, then Definition \ref{bb} coincides with the definition of  balanced pairs of triangulated categories (cf. \cite{FHZZ}). If the categoriy $\mathcal{C}$ are abelian, then Definition \ref{bb} coincides with the definition of  balanced pairs of abelian categories (cf. \cite{C}).

(2) By Remark \ref{xxx}, the condition (B1) may be rephrased as: any $\mathcal{X}$-resolution of $M$ is $C(-,\mathcal{Y})$-exact. Similarly, the condition (B2) may be rephrased as: any $\mathcal{Y}$-coresolution of $N$ is $C(\mathcal{X},-)$-exact.

\end{remark}
Our first main result is the following.
\begin{theorem}\label{main}
Let $\mathcal{C}$ be an extriangulated category with a negative first extension. The assignments
$$\Psi:(\mathcal{X},\mathcal{Y})\longmapsto\xi_{\mathcal{X}}=\xi^{\mathcal{Y}}~~~ and ~~~\Phi:\xi\longmapsto(\mathcal{P}(\xi),\mathcal{I}(\xi))$$
give mutually inverse bijections between the following classes:
\begin{itemize}
\rm\item[(1)] Balanced pairs $(\mathcal{X},\mathcal{Y})$ in $\mathcal{C}$.
 \rm\item[(2)]Proper classes $\xi$ in $\mathcal{C}$ with enough {$\xi$-projectives} and enough {$\xi$-injectives}.
 \end{itemize}
\begin{proof}

 Let $(\mathcal{X},\mathcal{Y})$ be a balanced pair in $\mathcal{C}$. Then $ \xi_{\mathcal{X}}=\xi^{\mathcal{Y}} $ is the desired proper class such that ${\mathcal{X}=\mathcal{P}(\xi_{\mathcal{X}}})$ and $ {\mathcal{Y}=\mathcal{I}(\xi^{\mathcal{Y}}})$ by Proposition \ref{aag}. Conversely, if ${\xi}$ is a proper class in ${\mathcal{C}}$ with enough {$\xi$-projectives} and enough {$\xi$-injectives}. We put $(\mathcal{X},\mathcal{Y})=(\mathcal{P}(\xi),\mathcal{I}(\xi))$ is a balanced pair by  Proposition \ref{aag}.

  For any balanced pair $(\mathcal{X},\mathcal{Y})$, it is easy to check that $$\Phi\Psi(\mathcal{X},\mathcal{Y})=\Phi(\xi_{\mathcal{X}}=\xi^{\mathcal{Y}})=(\mathcal{P}(\xi_{\mathcal{X}}),\mathcal{I}(\xi^{\mathcal{Y}}))=(\mathcal{X},\mathcal{Y}).$$
On the other hand, if $\xi$ is a proper class in  $\mathcal{C}$ with enough {$\xi$-projectives} and enough {$\xi$-injectives}, it is easy to see that $$\Psi\Phi(\xi)=\Psi(\mathcal{P}(\xi),\mathcal{I}(\xi))=(\mathcal{P}(\xi_{\mathcal{X}}),\mathcal{I}(\xi^{\mathcal{Y}}))=\xi.$$  	
This completes the proof.
\end{proof}

\end{theorem}
By applying Theorem \ref{main} to triangulated categories, we get the following.
\begin{corollary}\rm{\cite[Theorem 1.1]{FHZZ}}
Let $\mathcal{C}$ as a triangulated category. The assignments $\Psi:(\mathcal{X},\mathcal{Y})\longmapsto\xi_{\mathcal{X}}=\xi^{\mathcal{Y}}$ and $\Phi:\xi\longmapsto(\mathcal{P}(\xi),\mathcal{I}(\xi))$
give mutually inverse bijections between the following classes:
\begin{itemize}
\rm\item[(1)] Balanced pairs $(\mathcal{X},\mathcal{Y})$ in $\mathcal{C}$.
 \rm\item[(2)]Proper classes $\xi$ in $\mathcal{C}$ with enough {$\xi$-projectives} and enough {$\xi$-injectives}.
 \end{itemize}
\begin{proof}
This follows from Theorem \ref{main} and Remark \ref{111} and  Remark \ref{222}.
\end{proof}
\end{corollary}

By applying Theorem \ref{main} to abelian categories, we get the following.
\begin{corollary}\rm{\cite[Theorem 2.2]{WLH}}
Let $\mathcal{C}$ as an abelian category, there exists a one-to-one correspondence between balanced pairs and Quillen exact structures $\xi$ with enough {$\xi$-pojectives} and enough {$\xi$-injectives}.
\begin{proof}
This follows from Theorem \ref{main} and Remark \ref{111} and Remark \ref{222}.
\end{proof}
\end{corollary}

\begin{corollary}
\rm Let $(\mathcal{X},\mathcal{Y})$ be a balanced pair in an extriangulated category $\mathcal{C}$ with a negative first extension. Then $\xi:=\xi_{\mathcal{X}}=\xi^{\mathcal{Y}}$
is a proper class in $\mathcal{C}$. With the notation above, we set $\mathbb{E}_\xi:=\mathbb{E}|_\xi$, $\mathbb{E}^{-1}_\xi:=\mathbb{E}^{-1}|_\xi$, that is,
$$\mathbb{E}_\xi(C, A)=\{\delta\in\mathbb{E}(C, A)~|~\delta~ \textrm{is realized as an $\mathbb{E}$-triangle}\xymatrix{A\ar[r]^x&B\ar[r]^y&C\ar@{-->}[r]^{\delta}&}~\textrm{in}~\xi\}$$ for any $A, C\in\mathcal{C}$, and $\mathfrak{s}_\xi:=\mathfrak{s}|_{\mathbb{E}_\xi}$.
Hence $(\mathcal{C},\mathbb{E}_\xi,\mathbb{E}^{-1}_\xi,\mathfrak{s}_\xi)$ is an extriangulated category with a negative first extension.
\begin{proof}
This follows from Theorem \ref{main} and Theorem 3.2 in \cite{HZZ}.
\end{proof}
\end{corollary}

\section{Glued balanced pairs }
We always assume that any extrianglated category satisfies the (WIC) condition,
see \cite[Condition 5.8]{NP}. We briefly recall of the concepts and basic properties of recollements of extriangulated categories from \cite{WWZ}.
We omit some details here, but the reader can find them in \cite{WWZ}.

\begin{definition}\label{recollement}{\rm \cite[Definition 3.1]{WWZ}}
Let $\mathcal{A}$, $\mathcal{B}$ and $\mathcal{C}$ be three extriangulated categories. A \emph{recollement} of $\mathcal{B}$ relative to
$\mathcal{A}$ and $\mathcal{C}$, denoted by ($\mathcal{A}$, $\mathcal{B}$, $\mathcal{C}$), is a diagram
\begin{equation}\label{recolle}
  \xymatrix{\mathcal{A}\ar[rr]|{i_{*}}&&\ar@/_1pc/[ll]|{i^{*}}\ar@/^1pc/[ll]|{i^{!}}\mathcal{B}
\ar[rr]|{j^{\ast}}&&\ar@/_1pc/[ll]|{j_{!}}\ar@/^1pc/[ll]|{j_{\ast}}\mathcal{C}}
\end{equation}
given by two exact functors $i_{*},j^{\ast}$, two right exact functors $i^{\ast}$, $j_!$ and two left exact functors $i^{!}$, $j_\ast$, which satisfies the following conditions:
\begin{itemize}
  \item [(R1)] $(i^{*}, i_{\ast}, i^{!})$ and $(j_!, j^\ast, j_\ast)$ are adjoint triples.
  \item [(R2)] $\Im i_{\ast}=\Ker j^{\ast}$.
  \item [(R3)] $i_\ast$, $j_!$ and $j_\ast$ are fully faithful.
  \item [(R4)] For each $X\in\mathcal{B}$, there exists a left exact $\mathbb{E}_\mathcal{B}$-triangle sequence
  \begin{equation}\label{first}
  \xymatrix{i_\ast i^! X\ar[r]^-{\theta_X}&X\ar[r]^-{\vartheta_X}&j_\ast j^\ast X\ar[r]&i_\ast A}
   \end{equation}
  with $A\in \mathcal{A}$, where $\theta_X$ and  $\vartheta_X$ are given by the adjunction morphisms.
  \item [(R5)] For each $X\in\mathcal{B}$, there exists a right exact $\mathbb{E}_\mathcal{B}$-triangle sequence
  \begin{equation}\label{second}
  \xymatrix{i_\ast\ar[r] A' &j_! j^\ast X\ar[r]^-{\upsilon_X}&X\ar[r]^-{\nu_X}&i_\ast i^\ast X&}
   \end{equation}
 with $A'\in \mathcal{A}$, where $\upsilon_X$ and $\nu_X$ are given by the adjunction morphisms.
\end{itemize}
\end{definition}

\begin{remark}\label{77}(1) If the categories $\mathcal{A}$, $\mathcal{B}$ and $\mathcal{C}$ are abelian, then Definition \ref{recollement} coincides with the definition of recollement of abelian categories (cf. \cite{MH,P,FP}).

(2) If the categories $\mathcal{A}$, $\mathcal{B}$ and $\mathcal{C}$ are triangulated, then Definition \ref{recollement} coincides with the definition of recollement of triangulated categories (cf. \cite{BBD}).

(3) There is an example of recollement of an extriangulated category which is neither abelian nor triangulated, for more details, see \cite{WWZ} and  \cite{HHZ}.

\end{remark}

We collect some properties of a recollement of extriangulated categories, which will be used in the sequel.
\begin{lemma}\label{CY}\rm{\cite[Proposition 3.3]{WWZ}} Let ($\mathcal{A}$, $\mathcal{B}$, $\mathcal{C}$) be a recollement of extriangulated categories as \rm{(\ref{recolle})}. Then

$(1)$ all the natural transformations
$$i^{\ast}i_{\ast}\Rightarrow\Id_{\A},~\Id_{\A}\Rightarrow i^{!}i_{\ast},~\Id_{\C}\Rightarrow j^{\ast}j_{!},~j^{\ast}j_{\ast}\Rightarrow\Id_{\C}$$
are natural isomorphisms.

$(2)$ $i^{\ast}j_!=0$ and $i^{!}j_\ast=0$.

$(3)$  if $i^{\ast}$ is exact, then $j_{!}$ is  exact.

$(3')$ if $i^{!}$ is exact, then $j_{\ast}$ is exact.

\end{lemma}
\begin{lemma}\label{rr}\rm Let $(\mathcal{A}, \mathcal{B},\mathcal{C})$ be a recollement of extriangulated categories.

\label{rr0}$(1)$ If ${\mathcal{X}}$ is a \emph{strongly exact-contravariantly finite} subcategory of ${\mathcal{B}}$ and ${j_{!}j^{\ast}}\mathcal{X}\subseteq\mathcal{X}$, then ${{j^{\ast}}\mathcal{X}}$ is a \emph{strongly exact-contravariantly finite} subcategory of ${\mathcal{C}}$ ;

$(2)$ If ${\mathcal{Y}}$ is a \emph{strongly exact-covariantly finite} subcategory of ${\mathcal{B}}$ and ${j_{\ast}j^{\ast}}\mathcal{Y}\subseteq\mathcal{Y}$, then ${{j^{\ast}}\mathcal{Y}}$ is a \emph{strongly exact-covariantly finite} subcategory of ${\mathcal{C}}$.

\begin{proof}
$(1)$ For any $C$ in ${\mathcal{C}}$, $j_{\ast}C\in{\mathcal{B}}$. There exists an $\mathbb{E}_\mathcal{B}$-triangle $\xymatrix{K\ar[r]^{f}&X_0\ar[r]^{g}&j_{\ast}C\ar@{-->}[r]^{}&}$ with $X_0\in\mathcal{X}$, such that the sequence of abelian group $$0\rightarrow\mathcal{C}(X,K)\rightarrow\mathcal{C}(X,X_0)\rightarrow\mathcal{C}(X,j_{\ast}C)\rightarrow 0$$ is exact in ${\rm Ab},$ for any $X\in\mathcal{X}$, since ${\mathcal{X}}$ is a \emph{strongly exact-contravariantly finite} subcategory of ${\mathcal{B}}$. Because $j^{\ast}$ is exact, note that $j^{\ast}j_{\ast}\Rightarrow\rm{\Id}_\mathcal{C}$ is a natural isomorphism, applying $j^{\ast}$ to the above $\mathbb{E}_\mathcal{B}$-triangle, we obtain an $\mathbb{E}_\mathcal{C}$-triangle $$\xymatrix{j^{\ast}K\ar[r]^{j^{\ast}f}&j^{\ast}X_0\ar[r]^{j^{\ast}g}&C\ar@{-->}[r]^{}&}$$
with $j^{\ast}X_0\in j^{\ast}\mathcal{X}$. We need to show that the sequence of abelian group
$$0\rightarrow\mathcal{C}(j^{\ast}X,j^{\ast}K)\rightarrow\mathcal{C}(j^{\ast}X,j^{\ast}X_0)\rightarrow\mathcal{C}(j^{\ast}X,C)\rightarrow 0$$ is exact in ${\rm Ab},$
for any $X\in\mathcal{X}$. Since $(j_!, j^\ast, j_\ast)$ is a adjoint triple, there exists a commutative diagram
$$\xymatrix@C=1.8cm{
\mathcal{C}(j_!j^{\ast}X,K)\ar[r]^{\mathcal{C}(j_!j^{\ast}X, f)}\ar[d]^{\cong} &\mathcal{C}(j_!j^{\ast}X,X_0)\ar[d]^{\cong} &\\
\mathcal{C}(j^{\ast}X,j^{\ast}K)\ar[r]^{\mathcal{C}(j^{\ast}X, j^{\ast}f)} &\mathcal{C}(j^{\ast}X,j^{\ast}X_0)\ar[r]^{\mathcal{C}(j^{\ast}X, j^{\ast}g)}\ar[d]^{\cong} & \mathcal{C}(j^{\ast}X,C) \ar[d]^{\cong}\\
 &\mathcal{C}(X,j_\ast j^{\ast}X_0)\ar[r]^{\mathcal{C}(X, j_\ast j^{\ast}g)}\ar[r] &  \mathcal{C}(X,j_\ast C).
}$$
Therefore, it suffices to show that the first row is monic and the third row is epic by the above commutative diagram. Since the $\mathbb{E}_\mathcal{B}$-triangle $$\xymatrix{K\ar[r]^{f}&X_0\ar[r]^{g}&j_{\ast}C\ar@{-->}[r]^{}&}$$ is
$\mathcal C(\mathcal{X},-)$-exact, ${j_{!}j^{\ast}}\mathcal{X}\subseteq\mathcal{X}$, then ${\mathcal{C}(j_!j^{\ast}X, f)}$ is monic.
On the other hand, we need to show that there exists a morphism $h:X\longrightarrow j_\ast j^{\ast}X_0$, such that $j_\ast j^{\ast}g\circ h=\gamma$ for any $\gamma:X\longrightarrow  j_\ast{C}$. Suppose that $\nu_{X_0}$ is the unit of the adjoint pair $(j^{\ast},j_\ast)$, we have the following commutative diagram
$$\xymatrix{
&&X  \ar@{-->}[ddll]_h \ar@{-->}[dll]_s \ar[d]^\gamma\\
X_0 \ar[d]_{\nu_{X_0}} \ar[rr]_{g} &&j_\ast{C} \ar@{=}[d]  &\\
j_\ast j^{\ast}X_0 \ar[rr]^{j_\ast j^{\ast}g }&& j_\ast{C}. &
}
$$
Since $\xymatrix{K\ar[r]^{f}&X_0\ar[r]^{g}&j_{\ast}C\ar@{-->}[r]^{}&}$ is
$\mathcal C(\mathcal{X},-)$-exact, then there exists a morphism $s\colon X\to X_0$ such that $gs=\gamma$. Therefore we just need to take $h={\nu_{X_0}}s$, that is, $\mathcal C({X},j_\ast j^{\ast}g )$ is epic.
This shows that ${{j^{\ast}}\mathcal{X}}$ is a \emph{strongly exact-contravariantly finite} subcategory of ${\mathcal{C}}$.

$(2)$ For any $C\in{\mathcal{C}}$ and $j_{!}C\in{\mathcal{B}}$, there exists an $\mathbb{E}_\mathcal{B}$-triangle $$\xymatrix{j_{!}C\ar[r]^{f}&Y_0\ar[r]^{g}&L\ar@{-->}[r]^{}&}$$ with $Y_0\in\mathcal{Y}$, such that the sequence of abelian group $$0\rightarrow\mathcal{C}(L,Y)\rightarrow\mathcal{C}(Y_0,Y)\rightarrow\mathcal{C}(j_{!}C,Y)\rightarrow 0$$ is exact in ${\rm Ab},$ for any $Y\in\mathcal{Y}$, since ${\mathcal{Y}}$ is a \emph{strongly exact-covariantly finite} subcategory of ${\mathcal{B}}$. Because $j^{\ast}$ is exact, note that $\rm{\Id}_\mathcal{C}\Rightarrow j^{\ast}j_{!}$ is a natural isomorphism, applying $j^{\ast}$ to the above $\mathbb{E}_\mathcal{B}$-triangle, we obtain an $\mathbb{E}_\mathcal{C}$-triangle
$$\xymatrix{ C\ar[r]^{j^{\ast}f}&j^{\ast}Y_0\ar[r]^{j^{\ast}g}&j^{\ast}L\ar@{-->}[r]^{}&}$$
with $j^{\ast}Y_0\in j^{\ast}\mathcal{Y}$. We need to show that the sequence of abelian group
$$0\rightarrow\mathcal{C}(j^{\ast}L,j^{\ast}Y)\rightarrow\mathcal{C}(j^{\ast}Y_0,j^{\ast}Y)\rightarrow\mathcal{C}(C,j^{\ast}Y)\rightarrow 0$$ is exact in ${\rm Ab},$
for any $Y\in\mathcal{Y}$. Since $(j_!, j^\ast, j_\ast)$ is a adjoint triple, there is a commutative diagram
$$\xymatrix@C=1.8cm{
\mathcal{C}(L,j_\ast j^\ast Y)\ar[r]^{\mathcal{C}(g, j_\ast j^\ast Y)}\ar[d]^{\cong} &\mathcal{C}(Y_0,j_\ast j^\ast Y)\ar[d]^{\cong} &\\
\mathcal{C}(j^{\ast}L,j^{\ast}Y)\ar[r]^{\mathcal{C}(j^{\ast}g, j^{\ast}Y)} &\mathcal{C}(j^{\ast}Y_0,j^{\ast}Y)\ar[r]^{\mathcal{C}(j^{\ast}f, j^{\ast}Y)}\ar[d]^{\cong} & \mathcal{C}(C,j^{\ast}Y) \ar[d]^{\cong}\\
 &\mathcal{C}(j_!j^\ast Y_0,Y)\ar[r]^{\mathcal{C}(j_!j^\ast f,Y)}\ar[r] &  \mathcal{C}(j_!C,Y).
}$$
Therefore, it suffices to show that the first row is monic and the third row is epic by the above commutative diagram. Since the $\mathbb{E}_\mathcal{B}$-triangle $\xymatrix{j_{!}C\ar[r]^{f}&Y_0\ar[r]^{g}&L\ar@{-->}[r]^{}&}$ is $\mathcal C(-,\mathcal{Y})$-exact, ${j_{\ast}j^{\ast}}\mathcal{Y}\subseteq\mathcal{Y}$, then ${\mathcal{C}(g, j_\ast j^\ast Y)}$ is monic.
On the other hand, we need to show that there exists a morphism $\beta:j_{!} j^\ast Y_0\longrightarrow Y$, such that $\beta\circ j_{!} j^\ast f=\delta$ for any $\delta:j_{!}C\longrightarrow Y$. Suppose that $\mu_{Y_0}$ is the counit of the adjoint pair $(j_!,j^{\ast})$, we have the following commutative diagram
$$\xymatrix{
Y&\\
j_!C \ar@{=}[d]\ar[u]^{\delta} \ar[rr]_{j_!j^{\ast}f} &&j_!j^{\ast}Y_0 \ar[d]^{\mu_{Y_0}}\ar@{-->}[ull]_{\beta}  &\\
j_!C\ar[rr]^{f }&&Y_0\ar@{-->}[uull]_{t}. &
}
$$
Since  $\xymatrix{j_{!}C\ar[r]^{f}&Y_0\ar[r]^{g}&L\ar@{-->}[r]^{}&}$ is $C(-,\mathcal{Y})$-exact, then there exists a morphism $t:Y_0\longrightarrow Y$, such that $tf=\delta$. Therefore we just need to take $\beta=t{\mu_{Y_0}}$, that is, $C(j_! j^{\ast}f, {Y})$ is epic.

This shows that ${{j^{\ast}}\mathcal{Y}}$ is a \emph{strongly exact-covariantly finite} subcategory of ${\mathcal{C}}$.
\end{proof}
\end{lemma}

\begin{lemma}\label{TTT}\rm Let $(\mathcal{A}, \mathcal{B},\mathcal{C})$ be a recollement of extriangulated categories.

$(1)$ If ${\mathcal{X}}$ is a \emph{strongly exact-contravariantly finite} subcategory of ${\mathcal{B}}$ and $\mathbb{E}_{\mathcal{A}}^{-1}(i^{\ast}\mathcal{X},\mathcal{A})=0$, ${i^{\ast}}$ is exact, then ${{i^{\ast}}\mathcal{X}}$ is a \emph{strongly exact-contravariantly finite} subcategory of ${\mathcal{A}}$ ;

$(2)$ If ${\mathcal{Y}}$ is a \emph{strongly exact-covariantly finite} subcategory of ${\mathcal{B}}$ and $\mathbb{E}_{\mathcal{A}}^{-1}(\mathcal{A},i^{!}\mathcal{Y})=0$, ${i^{!}}$ is exact, then ${{i^{!}}\mathcal{Y}}$ is a \emph{strongly exact-covariantly finite} subcategory of ${\mathcal{A}}$.

\begin{proof}
It is similar to the proof of Lemma \ref{rr}.
\end{proof}
\end{lemma}

Our second main result is the following.

\begin{theorem}\label{www}\rm Let $(\mathcal{A}, \mathcal{B},\mathcal{C})$ be a recollement of extriangulated categories. Suppose that $(\mathcal{X}, \mathcal{Y})$ is a balanced pair in $\mathcal{B}$ with ${j_{!}j^{\ast}}\mathcal{X}\subseteq\mathcal{X},~{j_{\ast}j^{\ast}}\mathcal{Y}\subseteq\mathcal{Y},~{i_{\ast}i^{\ast}}\mathcal{X}\subseteq\mathcal{X},~{i_{\ast}i^{!}}\mathcal{Y}\subseteq\mathcal{Y}~ .$

$(1)$ If $\mathbb{E}_{\mathcal{C}}^{-1}(\mathcal{C},j^{\ast}\mathcal{Y})=\mathbb{E}_{\mathcal{C}}^{-1}(j^{\ast}\mathcal{X},\mathcal{C})=0$, then $(j^{\ast}\mathcal{X}, j^{\ast}\mathcal{Y})$ is a balanced pair in $\mathcal{C}$;

$(2)$ If ${i^{\ast}}$ and ${i^{!}}$ are exact, $\mathbb{E}_{\mathcal{A}}^{-1}(i^{\ast}\mathcal{X},\mathcal{A})=\mathbb{E}_{\mathcal{A}}^{-1}(\mathcal{A},i^{!}\mathcal{Y})=0$,  then $(i^{\ast}\mathcal{X}, {i^{!}}\mathcal{Y})$ is a balanced pair in $\mathcal{A}$.

\begin{proof}
 We check the conditions of balanced pair.

$(1)$ (B1) This follows form Lemma \ref{rr}.

(B2) For any $C\in\mathcal{C}$, $j_{\ast}C\in\mathcal{B}$. Since $(\mathcal{X}, \mathcal{Y})$ is a balanced pair in $\mathcal{B}$, it follows that there exists a $\mathcal{X}$-resolution as follows
  \begin{equation}\xymatrix@C=0.2cm{
  \cdots \ar[rr]^{} \ar[dr]_{}&& X_2\ar[rr]^{}\ar[dr]_{}&& X_1 \ar[rr]^{}\ar[dr]_{}&& X_0 \ar[rr]^{} &&j_{\ast}C \\
                &K_2\ar[ru]^{}&&K_1\ar[ru]^{} && K_0 \ar[ru]^{}  &&}
  \end{equation}
where $\xymatrix{K_0\ar[r]^{}& X_0\ar[r]^{}&j_{\ast}C\ar@{-->}[r]^{}&}$ and $\xymatrix{K_i\ar[r]^{}& X_i\ar[r]^{}&K_{i-1}\ar@{-->}[r]^{}&}$ in $\xi_{\mathcal{X}}, ~i\geq1$. Moreover, the $\mathcal{X}$-resolution is $\mathcal C(-,\mathcal{Y})$-exact.

Since $j^{\ast}$ is exact, note that $j^{\ast}j_{\ast}\Rightarrow\Id_\mathcal{C}$ is a natural isomorphism. Applying $j^{\ast}$ to (4.4), we obtain a complex as follows
\begin{equation}\xymatrix@C=0.2cm{
  \cdots \ar[rr]^{} \ar[dr]_{}&&j^{\ast} X_2\ar[rr]^{}\ar[dr]_{}&& j^{\ast}X_1 \ar[rr]^{}\ar[dr]_{}&&j^{\ast} X_0 \ar[rr]^{} &&C \\
                &j^{\ast}K_2\ar[ru]^{}&&j^{\ast}K_1\ar[ru]^{} && j^{\ast}K_0 \ar[ru]^{}  &&}
  \end{equation}

We claim that (4.5) is $\mathcal C(j^{\ast}\mathcal{X},-)$-exact. In fact, for any $X\in\mathcal{X}$, we have the following commutative diagram
$$\xymatrix@C=1.2cm{
 &\mathcal{C}(j^{\ast}X,j^{\ast}K_i)\ar[r]^{}\ar[d]^{\cong} & \mathcal{C}(j^{\ast}X,j^{\ast}X_i) \ar[r] \ar[d]^{\cong}& \mathcal{C}(j^{\ast}X,j^{\ast}K_{i-1})\ar[d]^{\cong}& \\
0\ar[r]^{} &\mathcal{C}(j_!j^{\ast}X,K_i)\ar[r]^{} &  \mathcal{C}(j_!j^{\ast}X,X_i)  \ar[r]& \mathcal{C}(j_!j^{\ast}X,K_{i-1})\ar[r]^{}& 0.
}$$
Note that ${j_{!}j^{\ast}}\mathcal{X}\subseteq\mathcal{X}$, then the first row is exact since the second row is exact in the above diagram. Besides, the sequence $0\rightarrow\mathcal{C}(j^{\ast}X,j^{\ast}K_0)\rightarrow\mathcal{C}(j^{\ast}X,j^{\ast}X_0)\rightarrow\mathcal{C}(j^{\ast}X,C)\rightarrow 0$ is exact by Lemma \ref{rr0} (1).

Next, we need to show (4.5) is $\mathcal C(-,j^{\ast}\mathcal{Y})$-exact by Remark \ref{222} (2). For any $Y\in\mathcal{Y}$, we have the following commutative diagram
$$\xymatrix@C=1.2cm{
 &\mathcal{C}(j^{\ast}K_{i-1},j^{\ast}Y)\ar[r]^{}\ar[d]^{\cong} & \mathcal{C}(j^{\ast}X_i,j^{\ast}Y) \ar[r] \ar[d]^{\cong}& \mathcal{C}(j^{\ast}K_{i},j^{\ast}Y)\ar[d]^{\cong}& \\
0\ar[r]^{} &\mathcal{C}(K_{i-1},j_{\ast}j^{\ast}Y)\ar[r]^{} &  \mathcal{C}(X_i,j_{\ast}j^{\ast}Y)  \ar[r]& \mathcal{C}(K_{i},j_{\ast}j^{\ast}Y)\ar[r]^{}& 0.
}$$
Note that ${j_{\ast}j^{\ast}}\mathcal{Y}\subseteq\mathcal{Y}$, then the first row is exact since the second row is exact in the above diagram. Moreover, consider the following  commutative diagram
$$\xymatrix@C=1.2cm{
 &\mathbb{E}_{\mathcal{C}}^{-1}(j^{\ast}K_{0},j^{\ast}Y)\ar[r]^{}&\mathcal{C}(C,j^{\ast}Y)\ar[r]^{} & \mathcal{C}(j^{\ast}X_0,j^{\ast}Y) \ar[r] \ar[d]^{\cong}& \mathcal{C}(j^{\ast}K_{0},j^{\ast}Y)\ar[d]^{\cong}& \\
 &&&  \mathcal{C}(X_0,j_{\ast}j^{\ast}Y)  \ar[r]& \mathcal{C}(K_{0},j_{\ast}j^{\ast}Y)\ar[r]^{}& 0.
}$$

Since the second row is exact in the above diagram, note that $\mathbb{E}_{\mathcal{C}}^{-1}(\mathcal{C},j^{\ast}\mathcal{Y})=0$, then the sequence $0\rightarrow\mathcal{C}(C,j^{\ast}Y)\rightarrow\mathcal{C}(j^{\ast}X_0,j^{\ast}Y)\rightarrow \mathcal{C}(j^{\ast}K_{0},j^{\ast}Y)\rightarrow 0$ is exact.

(B3) It is a dual of the proof of (B2). This completes the proof of (1).

$(2)$ (B1) This follows form Lemma \ref{TTT}.

(B2) Since ${{i^{\ast}}\mathcal{X}}$ is a \emph{strongly exact-contravariantly finite} subcategory of ${\mathcal{A}}$, there exists a ${i^{\ast}}\mathcal{X}$-resolution as follows
  \begin{equation}\xymatrix@C=0.3cm{
  \cdots \ar[rr]^{} \ar[dr]_{}&& i^{\ast}X_2\ar[rr]^{}\ar[dr]_{}&& i^{\ast}X_1 \ar[rr]^{}\ar[dr]_{}&& i^{\ast}X_0 \ar[rr]^{} &&A \\
                &K_2\ar[ru]^{}&&K_1\ar[ru]^{} && K_0 \ar[ru]^{}  &&}
  \end{equation}
where $\xymatrix{K_0\ar[r]^{}&i^{\ast}X_0\ar[r]^{}&A\ar@{-->}[r]^{}&}$ and $\xymatrix{K_i\ar[r]^{}& i^{\ast}X_i\ar[r]^{}&K_{i-1}\ar@{-->}[r]^{}&}$ in $\xi_{i^{\ast}\mathcal{X}}, ~i\geq1$, for any $A\in{\mathcal{A}}$.
It suffices to show that the resolution (4.6) is $\mathcal C(-,i^{!}\mathcal{Y})$-exact. Since $i_{\ast}$ is exact, applying $i_{\ast}$ to (4.6), we obtain a complex as follows
 \begin{equation}\xymatrix@C=0.2cm{
  \cdots \ar[rr]^{} \ar[dr]_{}&& i_{\ast}i^{\ast}X_2\ar[rr]^{}\ar[dr]_{}&&i_{\ast} i^{\ast}X_1 \ar[rr]^{}\ar[dr]_{}&& i_{\ast}i^{\ast}X_0 \ar[rr]^{} &&i_{\ast}A \\
                &i_{\ast}K_2\ar[ru]^{}&&i_{\ast}K_1\ar[ru]^{} &&i_{\ast} K_0 \ar[ru]^{}  &&}
  \end{equation}
Note that ${i_{\ast}i^{\ast}}\mathcal{X}\subseteq\mathcal{X}$, We claim that (4.7) is a $\mathcal{X}$-resolution of $i_{\ast}A $. In fact, it is easy to see the resolution (4.7) is $\mathcal C(\mathcal{X},-)$-exact since the resolution (4.6) is $\mathcal C(i^{\ast}\mathcal{X},-)$-exact and $(i^{\ast},i_{\ast})$ is an adjoint pair.
Since  $(\mathcal{X}, \mathcal{Y})$ is a balanced pair in $\mathcal{B}$, it follows that(4.7) is $C(-,\mathcal{Y})$-exact by Remark \ref{222} (2). Moreover, since $(i_{\ast},i^{!})$ is an adjoint pair, there exits the following two commutative diagrams:
 \begin{equation}\xymatrix@C=1.2cm{
 &\mathcal{C}(A,i^{!}Y)\ar[r]^{}\ar[d]^{\cong} & \mathcal{C}(i^{\ast}X_0,i^{!}Y) \ar[r] \ar[d]^{\cong}& \mathcal{C}(K_{0},i^{!}Y)\ar[d]^{\cong}& \\
0\ar[r]^{} &\mathcal{C}(i_{\ast}A,Y)\ar[r]^{} &  \mathcal{C}(i_{\ast}i^{\ast}X_0,Y)  \ar[r]& \mathcal{C}(i_{\ast}K_{0},Y)\ar[r]^{}& 0.
}  \end{equation}
and
 \begin{equation}\xymatrix@C=1.2cm{
 &\mathcal{C}(K_{i-1},i^{!}Y)\ar[r]^{}\ar[d]^{\cong} & \mathcal{C}(i^{\ast}X_i,i^{!}Y) \ar[r] \ar[d]^{\cong}& \mathcal{C}(K_{i},i^{!}Y)\ar[d]^{\cong}& \\
0\ar[r]^{} &\mathcal{C}(i_{\ast}K_{i-1},Y)\ar[r]^{} &  \mathcal{C}(i_{\ast}i^{\ast}X_i,Y)  \ar[r]& \mathcal{C}(i_{\ast}K_{i},Y)\ar[r]^{}& 0.
} \end{equation}
where the second rows are exact of (4.8) and (4.9). Thus, the resolution (4.6) is $C(-,i^{!}\mathcal{Y})$-exact.

(B3) It is a dual of the proof of (B2). This completes the proof of (2).
\end{proof}
\end{theorem}

By applying Theorem \ref{www} to abelian categories, we have the following.
\begin{corollary}\label{abe}\rm{\cite[Proposition 2.3]{FHY}} \rm Let $(\mathcal{A}, \mathcal{B},\mathcal{C})$ be a recollement of abelian categories. If $(\mathcal{X}, \mathcal{Y})$ is a balanced pair in $\mathcal{B}$ with ${j_{!}j^{\ast}}\mathcal{X}\subseteq\mathcal{X},~{j_{\ast}j^{\ast}}\mathcal{Y}\subseteq\mathcal{Y},~{i_{\ast}i^{\ast}}\mathcal{X}\subseteq\mathcal{X},~{i_{\ast}i^{!}}\mathcal{Y}\subseteq\mathcal{Y}$, then

$(1)$ $(j^{\ast}\mathcal{X}, j^{\ast}\mathcal{Y})$ is a balanced pair in $\mathcal{C}$;

$(2)$ $(i^{\ast}\mathcal{X}, {i^{!}}\mathcal{Y})$ is a balanced pair in $\mathcal{A}$.

\begin{proof}
Since $\mathcal{C}(\mathcal{X},-)$ and $\mathcal{C}(-,\mathcal{Y})$ are left exact, this assumptions of Theorem \ref{www}  are not necessary in abelian categories. This follows from Theorem \ref{www} and Remark \ref{222}.
\end{proof}
\end{corollary}

\textbf{Jian He}\\
Department of Mathematics, Nanjing University, 210093 Nanjing, Jiangsu, P. R. China\\
E-mail: \textsf{jianhe30@163.com}\\[0.3cm]
\textbf{Panyue Zhou}\\
College of Mathematics, Hunan Institute of Science and Technology, 414006 Yueyang, Hunan, P. R. China.\\
E-mail: \textsf{panyuezhou@163.com}


\begin{thebibliography}{99}
\bibitem[AET]{AET} T. Adachi, H. Enomoto, M. Tsukamoto. Intervals of $s$-torsion pairs in extriangulated categories with negative first extensions, arXiv: 2013.09549v1, 2021.

\bibitem[B]{B} A. Beligiannis. Relative homological algebra and purity in triangulated categories. J. Algebra 227: 268--361, 2000.

\bibitem[BBD]{BBD} A. Beilinson, J. Bernstein, P. Deligne. Faisceaux pervers, in: Analysis and topology on singular spaces, {I}, Luminy, 1981, {Ast\'erisque}, {vol. 100}, Soc. Math. France,
  Paris, 5--171,  1982.

\bibitem[C]{C}  X. Chen. Homotopy equivalences induced by balanced pairs. J. Algebra 324: 2718--2731, 2010.

\bibitem[FHY]{FHY} X. Fu, Y. Hu, H. Yao. The resolution dimensions with respect to balanced pairs in the recollement of abelian categories. J. Korean Math. Soc. 56(4): 1031--1048, 2019.

\bibitem[FHZZ]{FHZZ} X. Fu, J. Hu, D. Zhang, H. Zhu. Balanced pairs on triangulated categories.
arXiv: 2109.00932, 2021.


\bibitem[FP]{FP} V. Franjou, T. Pirashvili. Comparison of abelian categories recollements. Doc. Math. 9: 41--56, 2004.

\bibitem[HYF]{HYF} Y. Hu, H. Yao, X. Fu. Tilting objects in triangulated categories. Comm. Algebra 48: 410--429, 2020.

\bibitem[HZZ]{HZZ} J. Hu, D. Zhang, P. Zhou. Proper classes and Gorensteinness in extriangulated categories. J. Algebra 551: 23--60, 2020.

\bibitem[HHZ]{HHZ} J. He, Y. Hu, P. Zhou. Torsion pairs and recollements of extriangulated categories. arXiv: 2014.04924v1, 2021.

\bibitem[MH]{MH} X. Ma, Z. Huang. Torsion pairs in recollements of abelian categories, Front. Math. China, 13(4): 875--892, 2018.

\bibitem[MV]{MV} R. MacPherson, K. Vilonen. Elementary construction of perverse sheaves. Invent. Math. 84(2):  403--435, 1986.

\bibitem[NP]{NP} H. Nakaoka, Y. Palu. Extriangulated categories, Hovey twin cotorsion pairs and model structures. Cah. Topol. G\'{e}om. Diff\'{e}r. Cat\'{e}g. 60(2): 117--193, 2019.

\bibitem[NP1]{NP1} H. Nakaoka, Y. Palu. External triangulation of the homotopy category of exact quasi-category. arXiv: 2004.02479, 2020.

\bibitem[P]{P} C. Psaroudakis. Homological theory of recollements of abelian categories. J. Algebra, 398: 63--110, 2014.

\bibitem[WLH]{WLH} J. Wang, H. Li, Z. Huang, Applications of exact structures in abelian categories. Publ. Math. Debrecen, 88: 269-286, 2016.

\bibitem[WWZ]{WWZ} L. Wang, J. Wei, H. Zhang, Recollements of extriangulated categories. arXiv: 2012.03258v1, 2020.

\bibitem[ZZ]{ZZ} P. Zhou, B. Zhu. Triangulated quotient categories revisited. J. Algebra 502: 196--232, 2018.





\end{thebibliography}
\end{document}